\documentclass[a4paper,12pt]{amsart}
\usepackage{amscd}
\usepackage{amsmath}
\usepackage{amssymb}

\theoremstyle{plain}
\newtheorem{theorem}{Theorem}[section]
\newtheorem{lemma}[theorem]{Lemma}
\newtheorem{cor}[theorem]{Corollary}
\newtheorem*{nolabel}{Theorem {\ref{result}}}

\theoremstyle{definition}
\newtheorem{definition}[theorem]{Definition}
\newtheorem{remark}[theorem]{Remark}

\newtheoremstyle{notation}{14pt}{14pt}
     {}
     {}
     {\bfseries}
     {.}
     {\newline}
     {}

\theoremstyle{notation}
\newtheorem*{notation}{Notations}

\input xy
\usepackage[all]{xy}

\begin{document}

\def\pr{\mathbb{P}}
\def\com{\mathbb{C}}
\def\zet{\mathbb{Z}}
\def\real{\mathbb{R}}
\def\field{\mathbb{K}}
\def\razio{\mathbb{Q}}
\def\nat{\mathbb{N}}
\def\derxcat{{\bf D}(X)}
\def\dermcat{{\bf D}(M)}
\def\derccat{{\bf D}(C)}
\def\rhom{R{\mathrm{Hom}}}
\def\gistr{\mathbb{G}}
\def\qiso{\simeq_{\mathrm{q.iso}}}

\title{A semiorthogonal decomposition for Brauer Severi schemes}

\author{Marcello Bernardara}

\address{Laboratoire J. A. Dieudonn\'e, Universit\'e
de Nice Sophia Antipolis, Parc Valrore, 06108 Nice Cedex 2 (France).}

\email{bernamar\char`\@math.unice.fr}

\begin{abstract}
A semiorthogonal decomposition for the bounded derived category (the category of perfect complexes
in a non smooth case) of coherent sheaves on a Brauer Severi scheme is given. It relies on
bounded derived categories (categories of perfect complexes in a non smooth case) of suitably
twisted coherent sheaves on the base.
\end{abstract}
\maketitle

\section{Introduction}

In this paper we give a semiorthogonal decomposition of the bounded derived category (the category of
perfect complexes in a non smooth case) of coherent sheaves on a Brauer Severi scheme $f: X \to S$.
Very roughly, Brauer Severi schemes could be seen as a kind of twisted projective bundles.
This leads us to generalize the semiorthogonal decomposition given in \cite{orlov} for
projective bundles, by considering twisted sheaves on the base $S$ instead of untwisted ones.\par
Let us recall what happens in the case of a projective bundle. Let $S$ be a smooth projective
variety, $E$ a vector bundle of rank $r+1$ over $S$. We consider its projectivization
$p:X={\pr}(E) \to S$. We then have the following semiorthogonal decomposition for the bounded derived category $\derxcat$ of coherent sheaves on $X$.
\begin{theorem} {\bf{(Orlov)}}
Let ${\bf{D}}(S)_k$ be the full and faithful subcategory of $\derxcat$ whose objects are all objects
of the form $p^{\ast}A \otimes {\mathcal O}_X (k)$ for an object $A$ of ${\bf{D}}(S)$. Then the set of
admissible subcategories
$$({\bf{D}}(S)_0, \ldots, {\bf{D}}(S)_r)$$
is a semiorthogonal decomposition for the bounded derived category $\derxcat$ of coherent
sheaves on $X$.
\end{theorem}
\begin{proof}
This is \cite{orlov}, Theorem 2.6.
\end{proof}
The aim of the paper is to give the following generalization. Let $f: X \to S$ be a Brauer
Severi scheme of relative dimension $r$ over a locally notherian scheme $S$. Let $\alpha$
be the corresponding class
in $H^2 (S, {\gistr}_m)$. Let us
denote by $\derxcat$ the category of perfect complexes of coherent sheaves on $X$ and by
${\bf{D}}(S, \alpha)$ the category of perfect complexes of $\alpha$-twisted coherent sheaves on $S$.
Notice that in the smooth case they actually correspond to the bounded derived categories.
\begin{nolabel}
There exist admissible full subcategories ${\bf{D}}(S,X)_k$ of $\derxcat$, such that
${\bf{D}}(S,X)_k$ is equivalent to the category ${\bf{D}}(S, \alpha^{-k})$ for all $k$ in $\zet$.\par
The set of admissible subcategories
$$({\bf{D}}(S,X)_0, \ldots, {\bf{D}}(S,X)_r)$$
is a semiorthogonal decomposition for the category $\derxcat$ of perfect complexes of coherent
sheaves on $X$.
\end{nolabel}
It will be clear in the proof of the theorem that the construction of the full admissible
subcategories ${\bf{D}}(S,X)_k$ is strictly related to the definition of the full
admissible sucategories ${\bf{D}}(S)_k$ in Orlov's theorem.\par
The paper is organized as follows: in section \ref{twistdefine} we give the definition of twisted
sheaves and we state the basic facts about their connection with Brauer Severi schemes. In section
\ref{twistderive} we recall basic facts about derived categories, categories of perfect
complexes and derived functors in twisted case, following \cite{caldararu}. In section
\ref{semiadmiss} we recall the definition of admissible subcategories and semiorthogonal decomposition
in a triangulated category and we state some basic results about it. The main Theorem and
its proof are given in section \ref{maintheorem}, together with a simple example.
\begin{notation}
All schemes considered are locally noetherian.\par
$S_{et}$ denotes the \'etale site of a scheme $S$.
For the definition of a (\'etale) site, see \cite{sga4mezzo} or \cite{milne}.\par
Given $U \to S$ in $Cov(S_{et})$, we denote $U''$ the fibered product $U \times_S U$
and $U'''$ the fibered product $U \times_S U \times_S U$. We call $p_1$ and $p_2$
the projections from $U''$ to $U$ and $q_{i,j}$ the
projections from $U'''$ to $U''$.\par
$f:X \to S$ is a Brauer Severi scheme of relative dimension $r$, that means that
$f$ is smooth and each fiber of $f$ is isomorphic to ${\pr}^r$.\par
$X_U$ denotes $f^{-1}(U)$ for $U$ in $Cov(S_{et})$. We use notations $X_{U}^{\prime \prime}$,
$X_{U}^{\prime \prime \prime}$, $p_{i,X}$ and $q_{i,j,X}$ in the natural way. Notice that
$X_{U}^{\prime \prime} = X_{U''}$.\par
Given a scheme $S$ we denote ${\bf{D}}(S)$ (respectively ${\bf{D}}(S, \alpha)$) the
triangulated category of perfect complexes of ($\alpha$-twisted) sheaves on $S$,
where $\alpha$ is an element of the Brauer group Br$(S)$.
\end{notation}

\section{Twisted sheaves}\label{twistdefine}

In this section, we give the definition of twisted sheaves and we
state the relationship between them and Brauer Severi schemes. We are working in \'etale topology,
but all can be defined and stated in analytic topology as well (see \cite{caldararu}, I, 1).
\begin{definition}
Let $S$ be a scheme with \'etale topology, $U \to S$ in $Cov(S_{et})$,
$\alpha \in \Gamma (U''',
\gistr_m)$ a 2-cocycle.\par
An $\alpha$-\it{twisted} sheaf \rm on $S$ is given by a
sheaf $E$ over $U$ and an
isomorphism $\phi: p_1^{\ast} E \cong p_2^{\ast}E$, such that
$$(q_{2,3}^{\ast}\phi) \circ (q_{1,2}^{\ast}\phi) = \alpha
(q_{1,3}^{\ast}\phi)$$
We say that such a sheaf is \it coherent \rm if $E$ is a coherent sheaf on $U$,
and we denote ${\mathrm{Mod}}(S, \alpha)$ the category of $\alpha$-twisted sheaves on $S$, ${\mathrm{Coh}}(S, \alpha)$ the category of coherent $\alpha$-twisted
sheaves on $S$ and ${\bf{D}}(S, \alpha)$ the category of perfect complexes of such sheaves.
\end{definition}
The category ${\mathrm{Mod}}(S, \alpha)$ does not change neither by refining the open cover
$U \to X$, nor by changing $\alpha$ by a cochain.
\begin{lemma}\label{modifybycochain}
If $\alpha$ and $\alpha'$ represent the same element of $H^2(S, \gistr_m)$, the categories
${\mathrm{Mod}}(S, \alpha)$ and ${\mathrm{Mod}}(S, \alpha')$ are equivalent.
\end{lemma}
\begin{proof}
This is \cite{caldararu}, Lemma 1.2.8. Indeed if $\alpha$ and $\alpha'$ are in the same
cohomology class they differ by a 1-cochain: $\alpha = \alpha' + \delta \gamma$. But then
sending any $\alpha'$-twisted sheaf $(E, \phi)$ to the $\alpha$-twisted sheaf $(E, \gamma
\phi)$ gives the required equivalence.
\end{proof}
\begin{remark}
Notice that in general the choice of the cochain $\gamma$ matters: different choices give
different equivalences. Since we are just interested in the existence of such equivalences
and not in a special one, in what follows this choice will not matter.
\end{remark}
Now we can see how twisted sheaves arise naturally when we
consider Brauer Severi schemes. Let $f: X \to S$ be a smooth morphism
between schemes such that each fibre is isomorphic to
${\pr}^r$. Then we call $X$ a Brauer Severi scheme of relative dimension
$r$ over $S$.\par
We can find a
covering $U \to S$ in $Cov(S_{et})$, such that $X_U = f^{-1} (U)$ is a projective bundle
over $U$ and $X_U \to X$ is a covering in $Cov(X_{et})$.
Then we have a local picture ${\pr}(E_U) \to U$, where $E_U$
is a locally free sheaf of rank $r+1$ on $U$ and we have an
isomorphism $\rho: {\pr}(E_U) \tilde{\to} X_U$. This fact is a classical
application of descent theory (\cite{dixexpo}, I, 8).\par
Consider the cartesian diagram
$$\xymatrix{
X_{U}^{\prime \prime \prime} \ar@<.7ex>[r] \ar[r] \ar@<-.7ex>[r] &
X_{U}^{\prime \prime} \ar@<.5ex>[r] \ar@<-.5ex>[r] & X_U
}$$
and call the projections $p_{i,X}$ and $q_{i,j,X}$. We have an isomorphism
$$\psi:= p_{1,X}^{\ast} \rho^{-1} \circ p_{2,X}^{\ast} \rho:
{\pr}(p_1^{\ast} E_U) \tilde{\longrightarrow} {\pr}(p_2^{\ast} E_U).$$
We would like to lift it to an isomorphism $\phi:p_1^{\ast} E_U \tilde{\to}
p_2^{\ast} E_U$.\par
Consider $U$ such that $p_1^{\ast} E_U$ and $p_2^{\ast} E_U$ can be trivialized. This
implies that $\psi$ is an automorphism of $U'' \times {\pr}^r$ and then
it gives a section of $PGL(r+1, U'')$. We can again refine $U$ in order to
obtain from it a section of $GL(r+1, U'')$, which will give us the required
isomorphism $\phi: p_1^{\ast}E_U \tilde{\to} p_2^{\ast}E_U$. Notice that this is not canonical
since it can be done up to a choice of an element of $\Gamma(U'', \gistr_m)$.\par
For this reason, we
have $(q_{1,2}^{\ast} \phi) \circ (q_{2,3}^{\ast} \phi) = \alpha_U (q_{1,3}^{\ast} \phi)$,
where $\alpha_U \in \Gamma(U''', \gistr_m)$.
We can see that $\alpha_U$ gives a cocycle and then $(E_U, \phi)$ is
an $\alpha$-twisted sheaf.\par
From now on, given a Brauer Severi scheme $f:X \to S$, we will consider
the $\alpha$-twisted sheaf $(E_U, \phi)$ described above and the category
${\bf{D}}(S, \alpha)$.
Notice that the choice of $\alpha_U$ could be modified by a 1-cochain, but, by
Lemma \ref{modifybycochain}, this would give an equivalent category. In fact
everything depends just on the cohomology class of $\alpha$.\par
The class $\alpha$ represents the obstruction to $f: X \to S$ to be a projective bundle.
To express this via cohomology, recall the exact sequence of sheaves over $S$:
$$1 \longrightarrow \gistr_m \longrightarrow GL(r+1) \longrightarrow
PGL(r+1) \longrightarrow 1.$$
It gives a long cohomology sequence:
$$\dots \to H^1(S, GL(r+1))
\to H^1(S, PGL(r+1)) \stackrel{\delta}{\to} H^2(S, {\gistr}_m)
\to \dots$$
and especially a connecting homomorphism
$\delta$.\par
Let $[X]$ be the cohomology class
of $X$ in $H^1 (S, PGL(r+1))$ and
$\alpha' := \delta ([X])$ in $H^2 (S, \gistr_m)$. If $\alpha'
= 0$, the class $[X]$ would lift to an element of $H^1 (S,GL(r+1))$, that is a rank $r+1$ vector
bundle on $S$. Since $X$ is not a projective bundle, $\alpha'$ is a nonzero element of
the cohomological Brauer group Br$'(S) := H^2(S, {\gistr}_m)$
and it is exactly the cohomology class $\alpha$ of the $\alpha_U$ described above.\par
As a projective bundle ${\pr}(E_U)$ over $U$, on $X_U$ there exists a tautological line
bundle ${\mathcal O}_{X_U}(1)$. We will also write
${\mathcal O}_{X_U}(k)$ for $k \in {\zet}$, clearly meaning ${\mathcal O}_{X_U}(-1) =
{\mathcal O}_{X_U}(1)^{\vee}$ and so on.\par Notice that the choice of the
bundle ${\mathcal O}_{X_U}(1)$ over $X_U$ depends on the choice of
$E_U$, moreover ${\mathcal O}_{X_U}(1)$ does not glue as a global untwisted sheaf
${\mathcal O}_X (1)$ on $X$.
However, the existence of a section for the morphism $f$ ensures the existence of
a global ${\mathcal O}_X(1)$.
\begin{lemma}\label{sectionbundle}
Let $f:X \to S$ be a Brauer Severi scheme. If $s:S \to X$ is a section of $f$, then
there exists a vector bundle $G$ on $S$ such that ${\pr}(G) \cong X \to S$.
\end{lemma}
\begin{proof}
The result is known, but since it is hard to find a reference, we give a proof.\par
Consider the diagram
$$\xymatrix{
X_{U}^{\prime \prime \prime} \ar@<.7ex>[rr]^{q_{i,j,X}} \ar@<-.7ex>[rr] \ar[rr] \ar@<.5ex>[d]^{f} &&
X_{U}^{\prime \prime} \ar@<.5ex>[rr]^{p_{1,X}} \ar@<-.5ex>[rr]_{p_{2,X}} \ar@<.5ex>[d]^{f} && X_U \ar@<.5ex>[d]^{f} \ar[r] & X \ar@<.5ex>[d]^{f}\\
U''' \ar@<.7ex>[rr] \ar@<-.7ex>[rr]_{q_{i,j}} \ar[rr] \ar@<.5ex>[u]^{s} && U'' \ar@<.5ex>[rr]^{p_1} \ar@<-.5ex>[rr]_{p_2} \ar@<.5ex>[u]^{s} && U \ar@<.5ex>[u]^{s} \ar[r] & S. \ar@<.5ex>[u]^{s}}
$$
Here $s$ and $f$ are improperly used to mean their pull-backs to $U$, $U''$ and $U'''$ in order to
keep a clearer notation.\par
We can choose ${\mathcal{O}}_{X_U}(1)$ such that $s^{\ast}{\mathcal{O}}_{X_U}(1) =
{\mathcal O}_U$.\par
Consider now $p_{1,X}^{\ast} {\mathcal{O}}_{X_U}(1)$ and $p_{2,X}^{\ast} {\mathcal{O}}_{X_U}(1)$,
the two pull-backs of ${\mathcal{O}}_{X_U}(1)$ to $X_{U''}$. There exists an invertible sheaf $L$
on $S$ such that $p_{1,X}^{\ast}{\mathcal O}_{X_U}(1) \cong p_{2,X}^{\ast}{\mathcal O}_{X_U}(1)
\otimes f^{\ast} L$. Since
$$s^{\ast} p_{i,X}^{\ast} {\mathcal{O}}_{X_U}(1) = {\mathcal{O}}_{U''}$$
we have $L$ trivial. We choose an isomorphism
$$\phi:p_{1,X}^{\ast}{\mathcal O}_{X_U}(1) \longrightarrow p_{2,X}^{\ast}{\mathcal O}_{X_U}(1)$$
such that $s^{\ast} \phi = {\mathrm{Id}}_{{\mathcal O}_{U''}}$.\par
The isomorphism $\phi$ satisfies an untwisted cocycle condition. Indeed,
$$s^{\ast}((q_{1,2,X}^{\ast}\phi)\circ(q_{2,3,X}^{\ast}\phi)\circ(q_{1,3,X}^{\ast}\phi)^{-1}) =
{\mathrm{Id}}_{{\mathcal O}_{U'''}}.$$
This shows that ${\mathcal{O}}_{X_U}(1)$
gives a global untwisted sheaf ${\mathcal O}_X (1)$ and that means $X$ is a projective
bundle over $S$.
\end{proof}

\section{Derived categories and functors in twisted case}\label{twistderive}

In this section we show what happens to most common derived functors when we consider the
category of perfect complexes of twisted sheaves on a scheme. We will state theorems
we need for the rest of the paper. Proofs and a more satisfying
description can be found in \cite{caldararu}. It is in fact an adaptation to twisted case of the
results of \cite{hrd}.\par
\begin{remark}
In a non smooth case, the following theorems can not be stated if
we work in the bounded derived category of coherent sheaves. In
order to extend the ideas to a more general context, we will deal
with categories of perfect complexes of ($\alpha$-twisted) sheaves.
In the smooth case, they turn out to be the same as bounded
derived categories of ($\alpha$-twisted) coherent
sheaves, but keep in mind that in the non smooth case what  we call here
${\bf{D}}(S)$ (resp. ${\bf{D}}(S,\alpha)$) is not the bounded derived category of
($\alpha$-twisted) coherent sheaves on $S$ but just a full triangulated subcategory.\par
A complete treatement of perfect complexes on a site
is given in \cite{sga6}. Everything is defined in the very
general context of fibered categories, hence all definitions
fit for twisted sheaves.
\end{remark}
\begin{theorem}
Let $f: X \to S$ be a morphism between schemes, let $\alpha, \alpha'$ be in $H^2 (S, \gistr_m)$, and ${\mathcal{AB}}$
be the category of abelian groups.
Then the following derived functors are defined:
$$\begin{array}{rl}
\underline{R{\mathrm{Hom}}} &: {\bf{D}} (S, \alpha)^{\circ} \times {\bf{D}} (S, \alpha') \longrightarrow
{\bf{D}}(S, \alpha^{-1} \alpha')\\
R{\mathrm{Hom}} &: {\bf{D}} (S, \alpha)^{\circ} \times {\bf{D}} (S, \alpha) \longrightarrow
{\bf{D}}^b({\mathcal{AB}})\\
\bigotimes_S &: {\bf{D}} (S, \alpha) \times {\bf{D}} (S, \alpha') \longrightarrow
{\bf{D}}(S, \alpha \alpha')\\
Lf^{\ast} &: {\bf{D}} (S, \alpha) \longrightarrow {\bf{D}}(X, f^{\ast} \alpha)\end{array}$$
If $f: X \to S$ is a projective lci (locally complete intersection) morphism, then we can define:
$$
Rf_{\ast} : {\bf{D}}(X, f^{\ast} \alpha) \longrightarrow {\bf{D}} (S, \alpha).
$$
\end{theorem}

\begin{proof} \cite{caldararu}, Theorem 2.2.6. Remark that asking $f$ to be lci projective
is too restrictive, but sufficient in our case. See \cite{sga6}, III, 4 for details.
\end{proof}

\begin{theorem} {\bf{(Projection Formula).}} Let $f: X \to S$ be a projective lci morphism
between schemes, $\alpha, \alpha' \in H^2 (S, \gistr_m)$. Then there is a natural functorial
isomorphism
$$Rf_{\ast} (F) \otimes_{S} G \tilde{\longrightarrow} Rf_{\ast} (F \otimes_X Lf^{\ast}G)$$
for $F \in {\bf{D}}(X, f^{\ast}\alpha)$ and $G \in {\bf{D}}(S, \alpha')$.
\end{theorem}

\begin{proof} \cite{caldararu}, Theorem 2.3.5.\end{proof}

\begin{theorem} {\bf{(Adjoint property of $Rf_{\ast}$ and $Lf^{\ast}$).}} Let $f: X \to S$ be
a projective lci morphism between schemes, $\alpha \in H^2 (S, \gistr_m)$. Then we have
$$R{\mathrm{Hom}}(Lf^{\ast}F, G) \tilde{\longrightarrow} R{\mathrm{Hom}}(F, Rf_{\ast}G)$$
for $F \in {\bf{D}}(S, \alpha)$ and $G \in {\bf{D}}(X, f^{\ast}\alpha)$.
\end{theorem}
\begin{proof} \cite{caldararu}, Theorem 2.3.9.\end{proof}
\begin{theorem} {\bf{(Flat Base Change).}} Let $f: X \to S$ be a projective lci morphism between schemes,
$\alpha \in H^2 (S, \gistr_m)$. Let $u:S' \to S$ be a flat morphism, let
$X' = X \times_S S'$ and $v$, $g$ projections in the cartesian square:
$$\xymatrix{
X' \ar[d]_{g} \ar[r]^{v} & X \ar[d]^{f}\\
S' \ar[r]_{u} & S.}
$$
Then there is a natural functorial isomorphism:
$$u^{\ast} Rf_{\ast} F \tilde{\longrightarrow} Rg_{\ast} v^{\ast}F$$
for $F \in {\bf{D}}(X, f^{\ast}\alpha)$.
\end{theorem}
\begin{proof} \cite{caldararu}, Theorem 2.3.10.\end{proof}

\section{Semiorthogonal decompositions}\label{semiadmiss}

Let $k$ be a field and $\bf D$ a $k$-linear triangulated category.
\begin{definition}
A full triangulated subcategory ${\bf D}' \subset {\bf D}$ is \it
admissible \rm if the inclusion functor $i: {\bf D}' \to {\bf D}$ admits
a right adjoint.
\end{definition}
\begin{definition}
The \it orthogonal complement \rm ${\bf D}'^{\perp}$ of ${\bf D}'$ in ${\bf D}$ is
the full subcategory of all objects $A \in {\bf D}$ such that ${\mathrm{Hom}}(B,A) = 0$
for all $B \in {\bf D}'$.
\end{definition}
We remark firstly that the orthogonal complement of an admissible subcategory is a triangulated
subcategory.\par
It can be shown that a full triangulated subcategory ${\bf D}' \subset {\bf D}$ is admissible if
and only if for all object $A$ of ${\bf D}$, there exists a distinguished triangle $B \to A \to C$
where $B \in {\bf D}'$ and $C \in {\bf D}'^{\perp}$, see \cite{bondal}. We also have
the following Theorem.
\begin{theorem}\label{admissiblesemiorth}
Let ${\bf D}'$ be a full triangulated sucategory of a triangulated category $\bf D$. Then
${\bf D}'$ is admissible if and only if $\bf D$ is generated by ${\bf D}'$ and
${\bf D}'^{\perp}$.
\end{theorem}
\begin{proof}
\cite{bondalkap}, Proposition 1.5, or \cite{bondal}, Lemma 3.1.
\end{proof}
Admissible subcategories occur when we have a fully faithful exact functor $F: {\bf D}'
\to {\bf D}$ that admits a right adjoint. To be precise, this functor defines an equivalence between
${\bf D}'$ and an admissible subcategory of $\bf D$.
\begin{definition}
An ordered sequence of admissible triangulated subcategories $\sigma = ({\bf D}_1, \ldots,
{\bf D}_n)$ is \it semiorthogonal \rm if, for all $i > j$, one has
${\bf D}_j \subset {\bf D}_i^{\perp}$. If $\sigma$ generates the category $\bf D$,
we call it a \it semiorthogonal decomposition \rm of $\bf D$.
\end{definition}
\begin{lemma}\label{gettingadmissible}
Let $\sigma = ({\bf D}_1, \ldots,
{\bf D}_n)$ be a set of ordered full subcategories of $\bf{D}$ sucht that ${\bf D}_j \subset
{\bf D}_i^{\perp}$ for all $i > j$ and $\sigma$ generates $\bf{D}$. Then ${\bf{D}}_i$ is
admissible for $i=1, \ldots,n$, and then $\sigma$ is a
semiorthogonal decomposition of $\bf{D}$.
\end{lemma}
\begin{proof}
Consider ${\bf D}_n$ and ${\bf D}_n^{\perp}$: they generate the category ${\bf D}$
and then they are admissible. In general,
consider ${\bf D}_i$ and ${\bf D}_i^{\perp}$ for $1 \leq i < n$: they generate the
category ${\bf D}_{i+1}^{\perp}$ and then they are admissible.
\end{proof}
For further information about admissible subcategories and semiorthogonal decomposition, see \cite{bondal, bondalkap, bondaorlov}.

\section{The main Theorem}\label{maintheorem}

Let now $f:X \to S$ be a Brauer Severi scheme of relative dimension $r$ and $\alpha$ in
Br$(S)$ the element associated to it as explained in section \ref{twistdefine}.
This section is dedicated to the proof of the following Theorem.
\begin{theorem}\label{result}
There exist admissible full subcategories ${\bf{D}}(S,X)_k$ of $\derxcat$, such that
${\bf{D}}(S,X)_k$ is equivalent to the category ${\bf{D}}(S, \alpha^{-k})$ for all $k$ in $\zet$.\par
The set of admissible subcategories
$$\sigma = ({\bf{D}}(S,X)_0, \ldots, {\bf{D}}(S,X)_r)$$
is a semiorthogonal decomposition for the category $\derxcat$ of perfect complexes of coherent
sheaves on $X$.
\end{theorem}
Recall that there exist a rank $r+1$ locally free sheaf $E_U$ on $U$, such that
$X_U = {\pr}(E_U)$ and that $E_U$ gives an
$\alpha$-twisted sheaf on $S$. Moreover on $X_U$ we have a
tautological line bundle ${\mathcal O}_{X_U}(1)$.\par
We split this section in three parts: in the first one we define the full subcategories
${\bf{D}}(S,X)_k$ of $\derxcat$ and we show the equivalence between ${\bf{D}}(S,X)_k$ and
${\bf{D}}(S, \alpha^{-k})$; all this is inspired by a construction by Yoshioka \cite{yoshioka}.
In the second one we show that the sequence $\sigma$ is indeed a
semiorthogonal decomposition. In the third one we give a simple example.

\subsection{Construction of ${\bf{D}}(S,X)_k$}
\begin{definition}
We define ${\bf{D}}(S,X)_k$, for $k \in \zet$, to be the full subcategory of ${\bf{D}}(X)$ generated
by objects $A$ such that
\begin{equation}\label{cohsxk}
A_{\vert X_U} \qiso f^{\ast} A_U \otimes {\mathcal O}_{X_U} (k)
\end{equation}
where $A_U$ is an object in ${\bf{D}}(U)$.
\end{definition}

\begin{lemma}\label{lemmafk}
For all $k$ in ${\zet}$, there is a functor
$$f^{\ast}_k:{\bf{D}}(S, \alpha^{-k}) \longrightarrow {\bf{D}}(S,X)_k$$
given by the association
\begin{equation}\label{funtore}
A_{\vert U} \mapsto f^{\ast} A_{\vert U} \otimes {\mathcal O}_{X_U} (k).
\end{equation}
\end{lemma}

\begin{proof}
Firstly, $X_U$ is the projective bundle ${\pr}(E_U)$ over $U$. We then have on $X_U$ the surjective
morphism $f^{\ast} E_U \twoheadrightarrow {\mathcal O}_{X_U} (1)$. Given $F$ an
$\alpha^{-1}$-twisted sheaf on $S$, we have the surjective morphism:
$$f^{\ast}(F_U \otimes E_U) = f^{\ast}F_U \otimes f^{\ast}E_U \twoheadrightarrow f^{\ast} F_U \otimes
{\mathcal O}_{X_U} (1).$$
Since $F_U$ and $E_U$ give respectively an $\alpha^{-1}$-twisted and an $\alpha$-twisted sheaf on $S$,
their tensor product $F_U \otimes E_U$ gives an untwisted sheaf on $S$.
Hence we see that $f^{\ast}F_U \otimes {\mathcal O}_{X_U} (1)$ gives a quotient of an
untwisted sheaf and hence an untwisted sheaf on $X$. It is now clear that given an object
$A$ in ${\bf{D}}(S, \alpha^{-1})$, the object given locally by (\ref{funtore}) belongs
to ${\bf{D}}(S,X)_1$.\par
The proof is similar for any $k$ in $\zet$.
\end{proof}

\begin{theorem}
The functor $f^{\ast}_k$ defined in Lemma \ref{lemmafk} is an equivalence between the
category ${\bf{D}}(S,\alpha^{-k})$ and the category ${\bf{D}}(S,X)_{k}$.
\end{theorem}
\begin{proof}
Given $A$ in ${\bf{D}}(S,X)_1$, consider the association over $U$
$$A_{\vert X_U} \mapsto Rf_{\ast} (A_{\vert X_U} \otimes {\mathcal O}_{X_U} (-1)).$$
We show that it gives a functor $\Lambda$ from ${\bf{D}}(S,X)_1$ to ${\bf{D}}(S,
\alpha^{-1})$ and that is the quasi-inverse functor of $f_1^{\ast}$.\par
Firstly, since $A$ is in ${\bf{D}}(S,X)_1$, on $X_U$ we have $A_{\vert
X_U} = f^{\ast} A_U \otimes {\mathcal O}_{X_U} (1)$, with $A_U$
in ${\bf{D}}(U)$. Evaluating $\Lambda$ on $A_{\vert X_U}$ we get
$$Rf_{\ast} (A_{\vert X_U} \otimes {\mathcal O}_{X_U} (-1)) = Rf_{\ast} f^{\ast} A_U.$$
Now use projection formula:
$$Rf_{\ast}f^{\ast} A_U = Rf_{\ast}
{\mathcal O}_X \otimes A_U.$$
We have $R^i f_{\ast} {\mathcal O}_X = 0$ for $i > 0$ and $f_{\ast} {\mathcal O}_X
= {\mathcal O}_S$, and then
\begin{equation}\label{atrapoco}
Rf_{\ast}f^{\ast}A_U \qiso A_U.
\end{equation}
It follows that $\Lambda$ associates to $A_{\vert X_U}$ the object $A_U$ in ${\bf{D}}(U)$.\par
At a level of coherent sheaves, by the same reasoning used in Lemma \ref{lemmafk},
we have the surjective morphism
$$f^{\ast}(F_U \otimes E_U) \twoheadrightarrow f^{\ast} F_U \otimes
{\mathcal O}_{X_U} (1).$$
Since $E_U$ is an $\alpha$-twisted sheaf on $S$, we can give to $F_U$ the structure of
$\alpha^{-1}$-twisted sheaf over $S$. This shows that
$\Lambda$ is actually a functor from the subcategory ${\bf{D}}(S,X)_1$ to the
category ${\bf{D}}(S,\alpha^{-1})$.\par
It is now an evidence by (\ref{atrapoco}) that $\Lambda$ and $f_1^{\ast}$ are each other
quasi-inverse.\par
The proof for $k \in \zet$ is similar.
\end{proof}
We then have constructed full subcategories ${\bf{D}}(S,X)_k$ of $\derxcat$, each one equivalent
to a category of perfect complexes of suitably twisted sheaves on $S$.\par
Notice that we have $f^{\ast}_0 = Lf^{\ast} = f^{\ast}$ since $f$ is flat. In this case,
everything is untwisted and we can recover a Lemma by Orlov.
\begin{lemma}\label{orlovlemma}
The functor $f^{\ast}: {\bf{D}}(S) \to \derxcat$ is a full and faithful embedding.
\end{lemma}
\begin{proof} \cite{orlov}, Lemma 2.1.
\end{proof}
The full subcategory of $\derxcat$ which is the image of ${\bf{D}}(S)$ under the functor $f^{\ast}$ is
in fact the category ${\bf{D}}(S,X)_0$ defined earlier.\par

\subsection{$\sigma$ is a semiorthogonal decomposition}
\begin{lemma}
For any $A$ in ${\bf{D}}(S,X)_k$ and $B$ in ${\bf{D}}(S,X)_n$ we have
$R{\mathrm{Hom}}(A,B) = 0$ for $r \geq k-n > 0$.
\end{lemma}
\begin{proof}
We have locally $A_{\vert X_U} = f^{\ast} A_U \otimes {\mathcal O}_{X_U} (k)$
and $B_{\vert X_U} = f^{\ast} B_U \otimes {\mathcal O}_{X_U} (n)$.\par
We have:
$$\underline{R{\mathrm{Hom}}}(A_{\vert X_U},B_{\vert X_U}) = \underline{R{\mathrm{Hom}}}
(f^{\ast} A_U \otimes {\mathcal O}_{X_U} (k), f^{\ast} B_U \otimes {\mathcal O}_{X_U} (n)) =$$
$$= \underline{R{\mathrm{Hom}}} (f^{\ast} A_U, f^{\ast} B_U \otimes {\mathcal O}_{X_U} (n-k)).$$
We now use the adjoint property of $f^{\ast}$ and $Rf_{\ast}$:
$$\underline{R{\mathrm{Hom}}} (f^{\ast} A_U, f^{\ast} B_U \otimes {\mathcal O}_{X_U} (n-k)) =
\underline{R{\mathrm{Hom}}} (A_U, Rf_{\ast} (f^{\ast} B_U \otimes {\mathcal O}_{X_U} (n-k)).$$
Now by projection formula
$$Rf_{\ast} (f^{\ast} B_U \otimes {\mathcal O}_{X_U} (n-k)) =
B_U \otimes Rf_{\ast}({\mathcal O}_{X_U} (n-k)).$$
We have $Rf_{\ast}({\mathcal O}_{X_U} (n-k)) = 0$ for $-r \leq n-k < 0$ and hence the sheaves $\underline{R
{\mathrm{Hom}}}(A,B)$ are zero.\par
Using the local to global Ext spectral sequence, we get the proof.
\end{proof}
We thus have an ordered set $ \sigma = ({\bf{D}}(S,X)_0, \ldots,{\bf{D}}(S,X)_r)$ of
orthogonal subcategories of $\derxcat$. Last step towards the proof of Theorem \ref{result} is
to show that it generates the whole category.\par
Consider the fiber square over $S$:
$$\xymatrix{
P:= X \times_S X \ar[d]_{q} \ar[r]^(.7){p} & X \ar[d]^{g}\\
X \ar[r]_{f} & S.}
$$
The morphism $g$ corresponds to $f: X \to S$. We call $P$ the product $X \times_S X$.\par
Consider the diagonal embedding $\Delta: X \to P$. It is a section for the
projection morphism $p: P \to X$. By Lemma \ref{sectionbundle}, there exists a
vector bundle $G$ on $X$ such that $P \cong {\pr}(G) \to X$.\par
Consider now on $P$ the surjective morphism: $p^{\ast} G \twoheadrightarrow {\mathcal O}_P
\to 0$. We also have the Euler short exact sequence on $P$:
$$0 \longrightarrow \Omega_{P / X} (1) \longrightarrow p^{\ast} G
\longrightarrow {\mathcal O}_{P}(1) \to 0.$$
Combining the exact sequence and the surjective morphism, we get a section of ${\mathrm{Hom}}
(\Omega_{P / X} (1), {\mathcal O}_P)$ whose zero
locus is the diagonal $\Delta$ of $P$. Remark that $\Omega_{P / X} (1) =
p^{\ast} \Omega_{X / S} \otimes {\mathcal O}_P (1)$ and
$$\Lambda^{k} (p^{\ast} \Omega_{X / S} \otimes {\mathcal O}_P (1)) =
p^{\ast} \Omega_{X / S}^{k} \otimes {\mathcal O}_P (k).$$
We get a Koszul resolution:
$$
0 \to p^{\ast} \Omega_{X / S}^r \otimes {\mathcal O}_P (r) \to
\ldots \to p^{\ast} \Omega_{X / S} \otimes {\mathcal O}_P (1) \to
{\mathcal O}_P \to {\mathcal O}_{\Delta} \to 0.
$$
By this complex we deduce that ${\mathcal O}_{\Delta}$ belongs, as an
element of the category ${\bf{D}}(P)$, to the subcategory generated from
\begin{equation}\label{insiemegen}
\lbrace p^{\ast} \Omega_{X / S}^r \otimes {\mathcal O}_P (r),
\ldots, p^{\ast} \Omega_{X / S} \otimes {\mathcal O}_P (1),
{\mathcal O}_X \boxtimes {\mathcal O}_X \rbrace
\end{equation}
by exact triangles and shifting.\par
Given $A$ an element of $\derxcat$, we remark that $A = Rq_{\ast}(p^{\ast}A \otimes {\mathcal O}_{\Delta})$.
Since all involved functors (pull-back, direct image and tensor product) are exact functors, $A$ belongs to the subcategory of $\derxcat$ generated by
$$\{ Rq_{\ast}(p^{\ast}(A \otimes \Omega^r_{X / S}) \otimes {\mathcal O}_P(r)), \ldots,
Rq_{\ast}(p^{\ast}(A \otimes \Omega_{X / S}) \otimes {\mathcal O}_P(1)), Rq_{\ast}p^{\ast}A \}.$$
\begin{lemma}
The object $Rq_{\ast}(p^{\ast}(A \otimes \Omega^k_{X / S}) \otimes {\mathcal O}_P(k))$
in $\derxcat$ belongs to the subcategory ${\bf{D}}(S,X)_k$.
\end{lemma}
\begin{proof}
We look at it in a local situation. In this case $X_U$ is a projective bundle over $U$,
and we have
$$q^{\ast}{\mathcal O}_{X_U}(k) = {\mathcal O}_P (k)_{\vert X_{U''}}.$$
This leads us to write locally:
\begin{equation}\label{passofinale}
Rq_{\ast} (p^{\ast} (A \otimes \Omega^k_{X / S}) \otimes {\mathcal O}_P(k))_{\vert X_{U''}}=
Rq_{\ast} ((p^{\ast} (A \otimes \Omega^k_{X / S}))_{\vert X_U} \otimes q^{\ast} {\mathcal O}_{X_U} (k)) =
\end{equation}
$$
= Rq_{\ast} (p^{\ast} (A \otimes \Omega^k_{X / S}))_{\vert X_U} \otimes {\mathcal O}_{X_U}(k)
=f^{\ast} Rg_{\ast} ((A \otimes \Omega^k_{X / S})_{\vert X_U}) \otimes {\mathcal O}_{X_U}(k)
$$
where we used projection formula and flat base change in the last two equalities.
Then we have an object locally of the form finally given in (\ref{passofinale}), and then it is
an object in ${\bf{D}} (S,X)_k$.
\end{proof}
We have shown that all objects $A$ in $\derxcat$ belong to the subcategory generated by the orthogonal
sequence $\sigma$. This implies, by Lemma \ref{gettingadmissible}, that
the subcategories ${\bf{D}}(S,X)_k$ are admissible and
then $\sigma$ is in fact a semiorthogonal decomposition of $\derxcat$. This completes the proof
of Theorem \ref{result}.

\subsection{An example}
We finally treat the simplest example of a Brauer Severi
scheme. Let $K$ be a field and $X$ a Brauer Severi variety over
the scheme $Spec(K)$. In this case Theorem
\ref{result} gives a very explicit semiorthogonal decomposition of
the bounded derived category $\derxcat$ of coherent sheaves on
$X$ in terms of central simple algebras over $K$.\par
The cohomological Brauer group of $Spec(K)$ is indeed the Brauer group
${\mathrm{Br}}(K)$ of the field $K$. The elements of
${\mathrm{Br}}(K)$ are equivalence classes of central simple
algebras over $K$ and its composition law is tensor product. To
each $\alpha$ in ${\mathrm{Br}}(K)$ corresponds the choice of a central simple algebra
over $K$.\par
Given the $\alpha$ corresponding to the Brauer Severi variety $X$, an $\alpha^{-1}$-twisted
sheaf is then a module over a properly chosen central simple algebra $A$,
and it is coherent if it is finitely generated. The category ${\bf{D}}(Spec(K),
\alpha^{-1})$ is the bounded derived category of finitely generated
modules over the algebra $A$. Concerning the element
$\alpha^{-k}$ in ${\mathrm{Br}}(K)$, just remind that the composition
law is tensor product, to see that we can choose $A^{\otimes k}$
to represent it. The construction of
${\bf{D}}(Spec(K), \alpha^{-k})$ is then straightforward. We can
state the following Corollary of the Theorem \ref{result}.
\begin{cor}
Let $K$ be a field, $X$ a Brauer Severi variety over $Spec(K)$ of dimension
$r$. Let $\alpha$ be the class of $X$ in ${\mathrm{Br}}(K)$ and $A$ a central simple
algebra over $K$ representing $\alpha^{-1}$.\par
The bounded derived category $\derxcat$ of coherent sheaves on $X$ has a semiorthogonal
decomposition $\sigma = ({\bf{D}}(K,X)_0, \ldots,{\bf{D}}(K,X)_r)$,
where ${\bf{D}}(K,X)_i$ is equivalent to the bounded derived
category of finitely generated $A^{\otimes i}$-modules.
\end{cor}

\end{document}